\documentclass[12pt,a4paper]{amsart}
\usepackage{amssymb}
\usepackage[final]{showkeys}
\usepackage{microtype}
\usepackage{color}
\usepackage{graphicx}
\usepackage[hmargin=3.5cm,vmargin={5cm,5cm}]{geometry}

\newcommand{\ds}{\displaystyle}
\numberwithin{equation}{section}

\allowdisplaybreaks

\begin{document}

    \title[]{A generalization of the Lomnitz logarithmic \\ creep
    law via Hadamard fractional calculus}

    \author{Roberto Garra$^1$}
     \address{${}^1$Dipartimento di Scienze Statistiche, 
    University of  Rome  ``La Sapienza", Italy.}
    \email{roberto.garra@sbai.uniroma1.it}

    \author{Francesco Mainardi$^2$}
    	 \address{${}^2$Dipartimento di Fisica e Astronomia (DIFA), 
    	 University of  Bologna ``Alma Mater Studiorum", and INFN, Italy}
	 \email{francesco.mainardi@bo.infn.it  \  (Corresponding Author)}
    
    \author{Giorgio Spada$^3$}
   \address{${}^3$Dipartimento  di Scienze Pure e Applicate
             (DiSPeA),  University of Urbino ``Carlo Bo", Italy.} 

\email{giorgio.spada@gmail.com}    


    \date{\today}

    \begin{abstract}
    We present a new approach based on linear integro-differential operators with logarithmic kernel   related to the Hadamard fractional calculus in order to generalize,
by a parameter $\nu \in (0,1]$,     
    the logarithmic creep law known in rheology as Lomnitz law
(obtained for $\nu=1$).    
  We  derive the constitutive stress-strain relation of this generalized model in a form that couples memory effects and time-varying viscosity.
    Then, based on the hereditary theory of linear viscoelasticity,
     we also derive the corresponding relaxation function by solving numerically a Volterra integral equation of the second kind. 
    So doing we provide a full characterization of the new model both in creep and in relaxation representation, where   the slow varying functions of logarithmic type play a fundamental role as required in processes of ultra slow kinetics.
      \\     
    {\bf   
     Chaos, Solitons and Fractals (2017):
Special Issue on Future Directions in Fractional Calculus.
    Guest Editors: Mark M. Meerschaert, Bruce J. West, Yong Zhou.
    Published on line  25 March 2017;  
    http://dx.doi.org/10.1016/j.chaos.2017.03.032
}
  \end{abstract}

    \maketitle

    \section{Introduction}
    
%    \section{The Lomnitz logarithmic creep law}

    In 1956  Lomnitz \cite{Lo1} introduced a logarithmic creep law 
     to treat the creep behaviour of igneous rocks. 
     The original Lomnitz creep law 
provides the strain response $\epsilon(t)$ to a constant stress $\sigma(t)=\sigma_0$ for $t\ge 0$ in the form     
     \begin{equation}\label{lo}
     \epsilon(t)= \frac{\sigma_0}{E_0}\left[1+q\ln(1+t/\tau_0)\right],
     \quad t\ge 0\,,
     \end{equation}
     where $E_0$ is the shear modulus, $\tau_0>0$  is a characteristic  time during which the transition from 
     elastic to creep-type deformation occurs 
     and $q$ is a positive non-dimensional constant.
  The logarithmic creep law suggested by Lomnitz on the basis of an empirical reasoning has found many applications in other 
     %% from the early
      papers of the same author \cite{Lo2,Lo3}.
      % to more recent investigations in the field of time-dependent non-Newtonian%  rheology \cite{Holm,Main}.
\newpage

In 1958 Jeffreys \cite{Jeffreys_GJRAS58}
 proposed a  power law of creep, generalizing the Lomnitz logarithmic law to broaden  the geophysical applications to  fluid-like materials including igneous rocks.
This generalized law, however, can be applied also to solid-like  viscoelastic materials as shown in a revisited version in 1984 by Strick \cite{Strick_JGR84}.
More  recently, in 2012 Mainardi and Spada \cite{Main}  have provided a full characterization  of the  rheological properties related to this general model thus including the original Lomnitz  creep law.   

We recall that in the linear   theory of viscoelasticity, based on the hereditary theory by Volterra,  
a viscoelastic body  is characterized by two distinct but interrelated material functions,
causal in time (i.e. vanishing for $t<0$):
the creep compliance $J(t)$ (the strain response to a unit step of stress)
and the relaxation modulus  $G(t)$ (the stress response to a unit step of strain).
For more details see e.g. the treatises by 
Christensen \cite{Christensen_BOOK82},
Pipkin \cite{Pipkin_BOOK86} and 
Mainardi \cite{Mainardi_BOOK10}.

By taking $J(0^+)=J_0 >0$ so that $G(0^+)= G_0 =1/J_0$,
the body is assumed to exhibit a non vanishing  instantaneous response both in the  creep  and in the relaxation tests.
 As a consequence, we find it convenient to introduce two non-dimensional quantities $\psi(t)$ and $\phi(t)$ as follows 
\begin{equation}
J(t)= J_0[1+\psi(t)]\,, \quad G(t) = G_0\, \phi(t)\,,
\end{equation}  
where $\psi(t)$ is a non-negative  increasing function with $\psi(0) =0$ and  
$\phi(t)$ is a non-negative decreasing function with $\phi(0)=1$. 
Henceforth, $\psi(t)$ and $\phi(t)$ will
be referred to as dimensionless creep function and relaxation function,
respectively.
Viscoelastic bodies can be distinguished in solid-like 
and fluid-like  whether $J(+\infty)$  is finite or infinite so that 
$G(+\infty)= 1/J(+\infty)$ is non zero or zero, correspondingly.  

It is quite common in linear viscoelasticity to require the existence of positive retardation and relaxation spectra for the material functions $J(t)$ and $G(t)$,
 as pointed out  by  Gross in his 1953 monograph on the mathematical structure of the theories of viscoelasticity \cite{Gross_BOOK53}.
  This implies,  as  formerly proved in 1973 by Molinari \cite{Molinari_1973}  and revisited  in 2005 by Hanyga \cite{Hanyga_2005}, see also  Mainardi's book    \cite{Mainardi_BOOK10},
that $J(t)$ and $G(t)$ turn out to be Bernstein and Completely Monotonic functions, respectively.       
For  their mathematical properties  the interested reader is referred to  the excellent monograph by Schilling et al. \cite{sch}.

    As pointed out e.g. in  \cite{Mainardi_BOOK10},
    the relaxation modulus $G(t)$      
    can be derived from the corresponding creep compliance $J(t)$  through the Volterra integral equation of the second kind
    \begin{equation}
    G(t)= \frac{1}{J_0}-\frac{1}{J_0}\,
    \int_0^t \!\frac{dJ}{dt'}G(t-t')\, dt'\, ;
    \label{integral-equation}
    \end{equation}
    as a consequence,
    the dimensionless relaxation function $\phi(t)$ 
    obeys the Volterra integral equation
    \begin{equation}
    \label{integraleqphi}
    \phi(t)=1-\int_0^t\frac{d\psi}{dt}\phi(t-t')dt'\,.
    \end{equation}
    
    Mainardi and Spada in \cite{Main}  have shown, both analytically and numerically,  that the relaxation function corresponding to the Lomnitz creep law  decays in time as the slow varying function 
    % proportional to 
  $1/\ln t$.
 
    Quite  recently  Pandey and Holm
\cite{Holm}
  have  discussed the meaning of the empirical Lomnitz logarithmic law in the framework of time-dependent non-Newtonian rheology, 
  where the stress-strain relation is
     \begin{equation}
          \sigma(t)= \eta(t)\dot{\epsilon}(t)\,, \quad t\ge 0\,,
     \end{equation} 
     where  $\eta(t)>0$  represents  a time-dependent   viscosity 
     coefficient.
     %%%
     In particular, they have shown that the stress-strain equation leading to the Lomnitz law is 
     \begin{equation}\label{eq0}
    \frac{q}{E_0}\,  \sigma(t)=  
     \left(1+\frac{t}{\tau_0}\right) \dot{\epsilon}(t), \quad t\ge 0\,,
     \end{equation}
     so that  the time evolution of  viscosity is represented by the  differential operator
 \begin{equation}
 \label{eqO1} 
    \widehat{O}_1^t  := \left(1+\frac{t }{\tau_0}\right)\, \frac{d}{dt}\,.
    \end{equation}     
       
      The starting point of Pandey and Holm \cite{Holm} is
        a spring-dashpot viscoelastic model (of Maxwell type) with  viscosity varying linearly in time  that provides a relaxation function decaying in time to zero as a negative power law.
     We note that  the stress-strain relation \eqref{eq0} indeed yields the Lomnitz creep law \eqref{lo}  for constant stress $\sigma(t) =\sigma_0 $, but 
     the power law for the  relaxation function  is not compatible with its  
     ultra-slow decay  of logarithmic type derived  by
     by Mainardi and Spada in \cite{Main}
       in the framework of the linear theory of viscoelasticity as solution of  the Volterra integral equation \eqref{integral-equation}.

     In Geophysics  there exists another approach to derive the Lomnitz creep law:
   it is due to Scheidegger \cite{Scheidegger_1970a, Scheidegger_1970b},
   who in 1970 proposed a non-linear stress-strain relation that reads
     \begin{equation}
     \dot \sigma (t) = 2\eta\, { \ddot {\epsilon}}(t) + 
     \beta \, {\dot{\epsilon}}^2(t)\,,
          \end{equation}
          where $\eta$ is the (constant) viscosity and $\beta$ a creep factor.
   The integration   of this differential equation  for a step input of stress
    $\sigma(t)=\sigma_0$ for $t\ge0$ leads to the Lomnitz creep law \eqref{lo} provided that
     $     2\eta/\beta = q \sigma_0/E_0$ and $\tau_0$ is a suitable time constant of integration.  This non-linear approach, however,  even if justified by the author  for some effects related to energy dissipation in rocks, has not found
     a validation in the literature up to nowadays. Furthermore, no investigation  for the  relaxation of stress under constant strain has been considered.

       The above discussion  implies  that in order to justify the Lomnitz creep law we have (at least, to our  knowledge) three possible ways:
       the standard one based on the linear hereditary Volterra theory,  
        the approach with  non constant viscosity followed by Pandey and Holm \cite{Holm} and the non-linear approach proposed by Scheidegger.
       
In this paper, starting with the Pandey and Holm \cite{Holm} approach, we set up an iterative
operational method based on operator \eqref{eqO1} which leads to the generalized Lomnitz law             
         \begin{equation}
         \label{lo-gen}
   \epsilon(t)= \frac{\sigma_0}{E_0}\left[1+q\frac{\ln^\nu\left(1+\frac{t}{\tau_0}\right)}{\Gamma(1+\nu)}\right]\,,
   \quad 0<\nu \le 1\,, \quad t\ge 0\,.
 \end{equation}
 We will then derive the corresponding relaxation function of this law
  by solving the related Volterra integral equation. 
    %%%%
   An interesting outcome of our  analysis is that the 
   resulting rheological model considers both a time varying viscosity and the memory effects required by the  hereditary theory of linear viscoelasticity. 
   
     We first give some mathematical preliminaries in the next Section and then in Section 3 we explain the meaning of our approach.
    In Section4 we discuss the implications of our work and finally conclusions
     are drawn in Section 5.
     
     For readers' convenience we add two appendices.
     In Appendix A  we recall the essentials  of the Hadamard fractional calculus on which our operational approach is based. In Appendix B we outline the numerical method adopted to solve the Volterra integral equation satisfied by the relaxation function of our generalized model.

    \section{Integro-differential operators with logarithmic kernels \label{ope}}
    
   In the recent paper by Beghin, Garra and Macci 
    \cite{jap}, an integro-differential operator with logarithmic kernel has been introduced in the context of correlated fractional negative binomial processes
    in statistics. 
    Using their notation, the time-evolution operator 
    $\widehat{O}^t_\nu$  acting on a sufficiently well-behaved function 
    $f(t)$
    is defined as
\begin{equation}
\begin{array}{ll}
   \widehat{O}^t_\nu f(t) &:={\displaystyle
    \frac{1}{\Gamma(n-\nu)}} \,  \times \\
    & {\displaystyle
    \int_{\frac{1-a}{b}}^{t}	
  \!\!  
    \ln^{n-1-\nu}\left(\frac{a+bt}{a+b\tau}\right)
    \bigg[\left((\frac{a}{b}+\tau)\frac{d}{d\tau}\right)^n f(\tau)\bigg]\frac{b}{a+b\tau}\, d\tau,}
    \end{array}
    \label{defO}
    \end{equation}
    for $n-1 < \nu < n \in \mathbb{N}$, $0<a \le 1$ and $b> 0$.
    
    A relevant property of this operator is given by the following result 
    \begin{equation}\label{pr}
    \widehat{O}^{t}_\nu \ln^\beta (a+bt) = \frac{\Gamma(\beta+1)}{\Gamma(\beta+1-\nu)}\ln^{\beta-\nu}(a+bt)
    \end{equation}
    for $\nu \in (0,1)$ and $\beta>-1 \setminus{\{0\}}$,
      see \cite{jap}, pag. 1057 for the details.
    Moreover we have that $ \widehat{O}^{t}_\nu \ const. = 0$.
%%%
    We refer to the Appendix A for a short survey about fractional-type operators with logarithmic kernel, starting from the so-called Hadamard fractional calculus.

    In analogy with the classical theory of  fractional calculus
(see e.g. the monograph by Kilbas, Srivastava and Trujillo \cite{kilbas}),    
    we  introduce the integral operator with logarithmic kernel 
acting on a sufficiently well-behaved function $f(t)$ as    
    \begin{equation}
    \widehat{I}^t_\alpha f(t) := \displaystyle \frac{1}{\Gamma(\alpha)}
       \int_{\frac{1-a}{b}}^{t}	\!\! \ln^{\alpha-1}\left(\frac{a+bt}{a+b\tau}\right)
          f(\tau)\frac{b}{a+b\tau}\, d\tau, \quad \alpha>0,
    \end{equation}
    so that we recognize 
    \begin{equation}
     \widehat{O}^{t}_\nu f(t)= \widehat{I}^t_{n-\nu} \left[\left((\frac{a}{b}+t)\frac{d}{dt}\right)^n f(t)\right]\, ,\quad n-1<\nu<n, \; n \in \mathbb{N}.
    \end{equation}
    Therefore, recalling the definition \eqref{eqO1} for the operator  
    $\widehat{O}^t_1$ with $\tau_0=1$ 
    and  taking  $\nu \in (0,1)$,  $a = b= 1 $ in the above equation,
     we obtain  
    \begin{equation}
    \widehat{O}^{t}_\nu f(t)= \widehat{I}^t_{1-\nu}
    \left[\left((1+t)\frac{d}{dt}\right) f(t)\right]= \widehat{I}^t_{1-\nu} \widehat{O}^t_1\, f(t),  \quad 0<\nu<1\,.
    \end{equation}
    
    We also observe that, by using the property \eqref{pr}, it can be proved that the composed Mittag-Leffler function 
    \begin{equation}
    E_{\nu,1}(-\ln^\nu(1+t))= \sum_{k=0}^\infty\frac{(-\ln^\nu(1+t))^k}{\Gamma(\nu k+1)}
    \end{equation}
    is an eigenfunction of the operator $\widehat{O}^t_\nu$.
Once again we find  the   Mittag-Leffler functions, that are known to  play a central role in the theory of fractional differential equations, see  e.g.  the recent monograph by Gorenflo et al \cite{book}. 
%%%%%
As a matter of fact, the composition of the classical Mittag-Leffler function with the power of a logarithmic function decays as a power of the logarithmic function.
 This behaviour can be interesting 
in the framework the so-called {\it ultra-slow kinetics} that includes 
phenomena of strong anomalous relaxation and diffusion and related stochastic processes.
On this respect, the interested reader may refer  to the seminal   paper of 1983 
by Sinai \cite{Sinai_TBA83} 
and then, e.g.  to 
Chechkin et al. \cite{Chechkin-et-al_EPL03}, 
Metzler and Klafter \cite{Metzler-Klafter_JPhysics04},
Mainardi et al. \cite{Mainardi-et-al_JVC07,Mainardi-et-al_JVC08}  
  and more recently to Wen Chen et al. \cite{WenChen-et-al_FCAA16}.

    We can observe that, heuristically, the integro-differential operator 
    $\widehat{O}_\nu^t$ can be obtained by means of a 
     time-change $t\rightarrow \left(\frac{a}{b}+t\right)$ starting from the definition of the Caputo fractional derivative.
    We also note  that the operator $\widehat{O}^t_{\nu}$ can be considered as a sort of \textit{fractional} counterpart of the differential operator $\widehat{O}_1^t$ appearing in the stress-strain equation governing the original Lomnitz's law, see Eq. (\ref{eq0}). 
    Moreover, we remark that for $a= 0$ and $b=1$ 
    the operator $\widehat{O}^t_\nu$ coincides with the regularized Hadamard fractional derivative (see Appendix A, \cite{bal} and \cite{jap} for more details).
    
    \section{A generalization of the Lomnitz law}
    We here consider a generalization of the Lomnitz stress-strain relation
    \eqref{eq0}  where the operator 
    $\widehat{O}_1^t$  in \eqref{eqO1} is replaced  by the integro-differential operator $\widehat{O}_\nu^t$ 
    defined in \eqref{defO} by taking $a=b=1$.
    Note that we have taken for the sake of simplicity and without loss of generality $\tau_0 =1$.
    
%%%%%%    
     This means that we are considering a generalized rheology with memory effects and time-varying viscosity, based on the following stress-strain relation
    \begin{equation}
    \begin{array}{ll}
   \frac{q}{E_0}\,  \sigma(t)&= \widehat{O}^t_\nu \epsilon(t)\\
     & = \displaystyle{\frac{1}{\Gamma(1-\nu)}
    \int_{0}^{t}	\ln^{-\nu}\left(\frac{1+t}{1+\tau}\right)
         \bigg[\left((1+\tau)\frac{d}{d\tau}\right) \epsilon(\tau)\bigg]\frac{1}{1+\tau}d\tau}\\
         & =\displaystyle{\frac{1}{\Gamma(1-\nu)}\int_{0}^{t}	\ln^{-\nu}\left(\frac{1+t}{1+\tau}\right)
                 \bigg[\widehat{O}_1^t \epsilon(\tau)\bigg]\frac{1}{1+\tau}d\tau},
         \end{array}
    \end{equation}
    with $\nu \in (0,1)$, $t\geq 0$.
    
   By using property \eqref{pr}, we obtain that the solution of this generalized fractional equation, setting $\sigma(t) =\sigma_0$, is 
   \begin{equation}
   \epsilon(t)= \frac{\sigma_0}{E_0}\left[1+q\frac{\ln^\nu\left(1+t\right)}{\Gamma(1+\nu)}\right],\quad 0<\nu\le 1, \quad  t\ge 0,
   \end{equation}
   which  represents our generalized Lomnitz creep law.

   The corresponding non-dimensional creep function for this model is 
       \begin{equation}
       \psi_\nu (t)= q\,\frac{\ln^\nu\left(1+t\right)}{\Gamma(1+\nu)}, 
       \quad 0<\nu\le 1, \quad  t\ge 0,
       \end{equation}
    so that the rate of creep is
       \begin{equation}
       \dot{\psi}_\nu (t)= q\,\nu  \frac{\ln^{\nu-1}(1+t)}{\Gamma(1+\nu)}\frac{1}{1+t}, \quad 0<\nu\le 1, \quad  t\ge 0.
       \end{equation}
       Therefore,
      in view of  \eqref{integral-equation},  the dimensionless relaxation function $\phi_\nu (t)$
     corresponding to our   generalized Lomnitz creep law obeys the integral equation
       \begin{equation}
       \label{VIE}
       \phi_\nu(t)=1-\frac{q \ \nu}{\Gamma(1+\nu)}
       \int_0^t \! \frac{\ln^{\nu-1}(1+t')}{1+t'}\phi_\nu(t-t')\, dt',
        \quad 0<\nu\le 1, \quad  t\ge 0,
       \end{equation}
       that is a Volterra  integral equation of the second kind with logarithmic kernel 
       if $0<\nu<1$ as expected.  Of course for $\nu=1$ 
      the   case of the Lomnitz model   is recovered, already investigated 
      by Mainardi and Spada \cite{Main}.

Figures 1 and 2 show the evolution of the dimensionless creep function $\psi_\nu(t)$ and
of the relaxation function $\phi_\nu(t)$ corresponding to the generalised
Lomnitz model, for various values of parameter $\nu$. From Figure 1, we
observe that the effect of decreasing the value of $\nu$ is that of
increasing the initial steepness of the creep function and of decreasing
its subsequent rate of variation. This behaviour stems from the analytical
expression for the rate of creep given by Eq. (3.4), valid for any value of
time. From Figure 2, a similar sensitivity to the value of $\nu$ is
observed for the relaxation function, although no closed-form expression is
available to corroborate this finding. However, the results in Figures 1
and 2 are both in agreement with the asymptotic representations given
below.

\begin{figure}[h!]
\centering
\includegraphics[width=12.0cm]{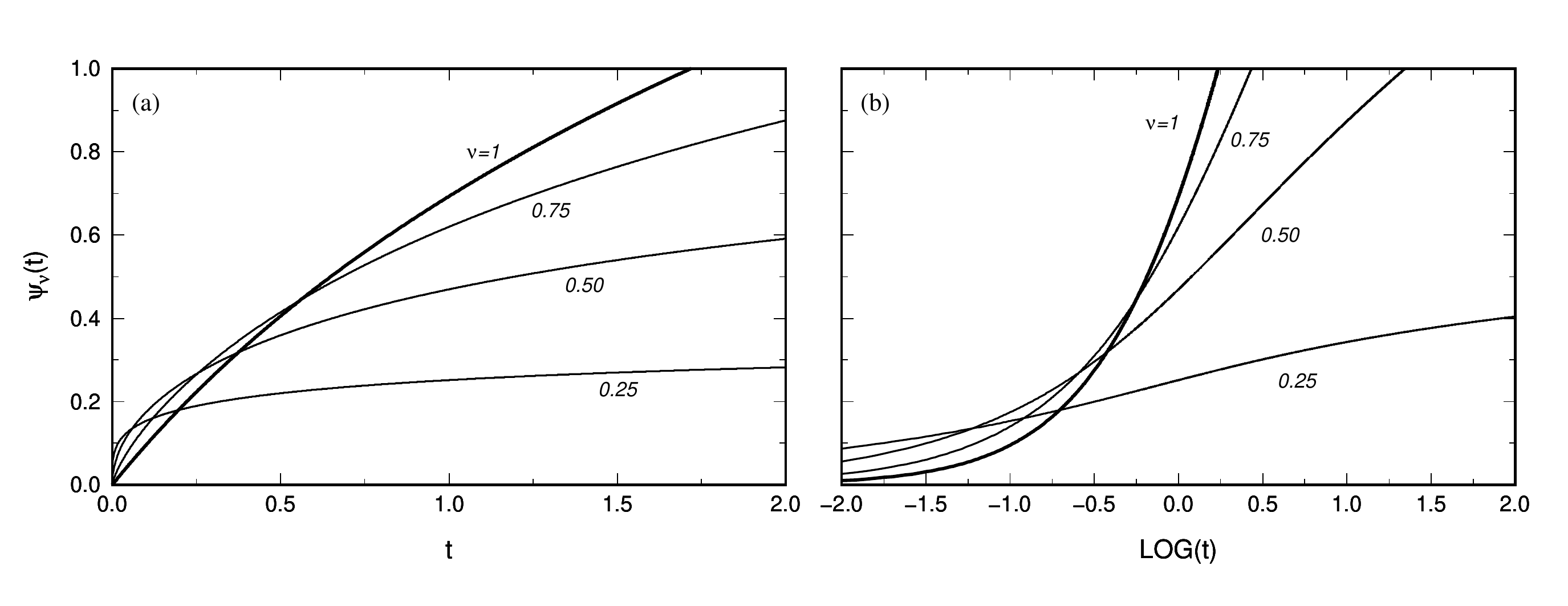}
\vskip-0.5truecm
\caption{The dimensionless creep function $\psi_\nu(t)$  of the generalized  Lomnitz model 
for $\nu = 0.25, 0.50, 0.75, 1 $ depicted versus time: 
(a) linear scale, (b) logarithmic scale.}
\end{figure}
 %%%%%%%%%
 \vskip-0.5truecm
 %%%%%%%%%%%%
   \begin{figure}[h!]
\centering
\includegraphics[width=12.0cm]{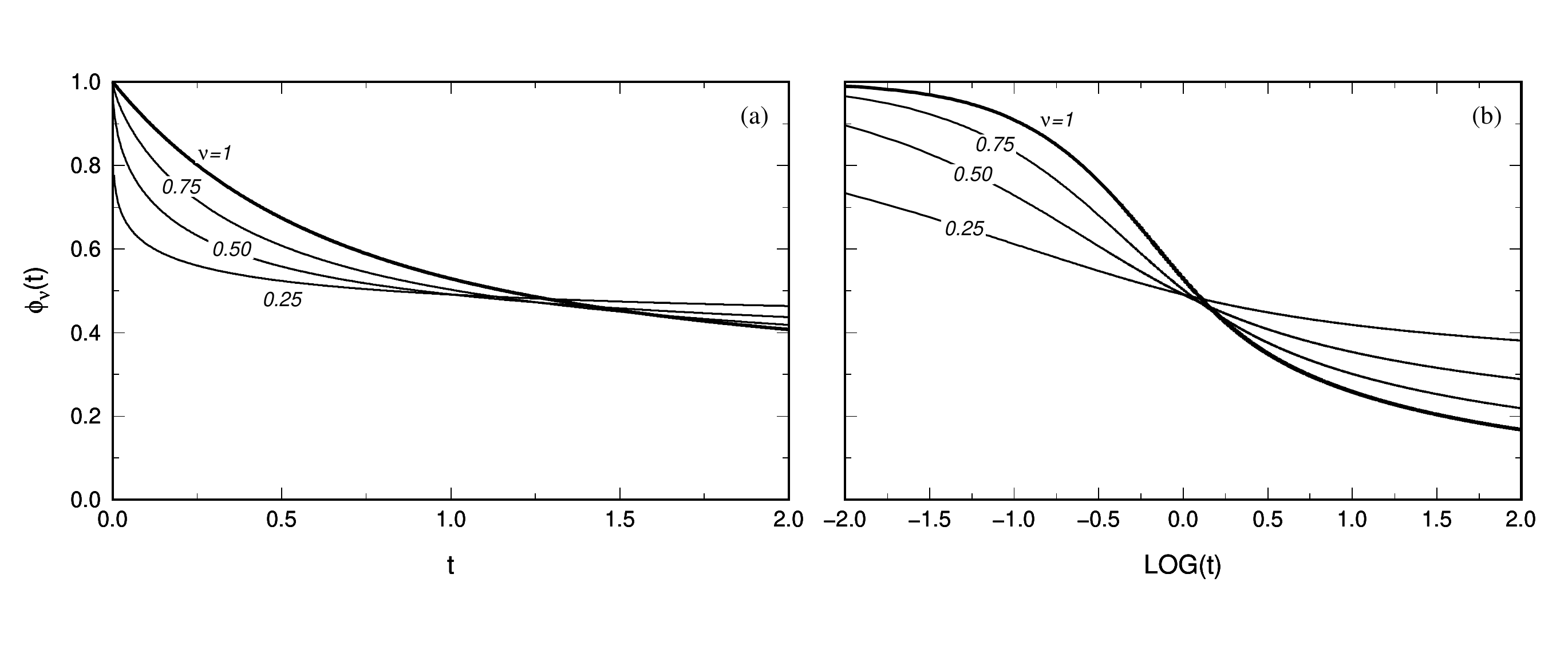}
\vskip-0.5truecm
\caption{The dimensionless relaxation function $\phi_\nu(t)$  of the generalized  Lomnitz model 
for $\nu = 0.25, 0.50, 0.75, 1 $ depicted versus time:
(a) linear scale, (b) logarithmic scale.}
\end{figure}
 
 Henceforth, we  report the asymptotic representations of the non-dimensional creep and relaxation functions taking $q=\tau_0=1$ for $0<\nu\le 1$
 for $t\to 0^+$ and $t\to +\infty$.
For $\psi_\nu(t)$  these representations are directly obtained from its analytical expression given by Eq. \eqref{lo-gen}.
 For
$\phi_\nu(t)$ the corresponding representations are properly derived from  the relation between the  Laplace transforms $\widetilde\phi_\nu(s)$ 
and $\widetilde\psi_\nu(s)$ that, by virtue of  the Volterra integral equation \eqref{VIE}, reads
\begin{equation}
 \widetilde\phi_\nu (s) = \frac{1}{s \,[1 +s \widetilde \psi_\nu(s)]}\,.
 \end{equation}
 Hence we get
 %%%%%%%%%% 
\begin{equation}
\psi_\nu(t) \sim
\left\{
\begin{array}{ll}
t^\nu/\Gamma(1+\nu)\,, & t\to 0^+\,, \\
\ln^\nu(t)/\Gamma(1+\nu)\,, & t\to +\infty\,,
\end{array}  
   \right.
   \end{equation}
   %%%%%%%%%%%%%%5
   and
\begin{equation}
\phi_\nu(t) \sim
\left\{
\begin{array}{ll}
1-t^{\nu}/\Gamma(1+\nu) \,, & t\to 0^+\,, \\
\Gamma(1+\nu)/\ln^\nu(t)\,, & t\to +\infty\,.
\end{array}  
   \right.
   \end{equation}   
  These asymptotic  representations have been checked 
  to fit with sufficiently good  approximation the numerical results in  suitable ranges of time.

  \section{Discussion}
  The applications of fractional rheologies in rock-physics is nowadays widely accepted, on the basis of empirical validations and rigorous theoretical models (see e.g. \cite{geo} and the references therein). Fractional models have the great advantage to take into account memory effects in relaxation processes, providing a more realistic picture of the physical mechanisms involved in rock-physics. The logarithmic Lomnitz creep law is essentially phenomenological, but a complete viscoelastic characterization has been recently provided in \cite{Main}. In view of the relevance of memory effects in relaxation processes, the advantage of our approach lies in a generalization of the Lomnitz law by considering an integro-differential counterpart of the stress-strain relation governing the original logarithmic creep law.
  This mathematical explorative model is based on the application of Hadamard-like fractional operators with logarithmic kernel. This choice is based on a physical reasoning: in this case the memory kernel should follow the logarithmic-type time of relaxation. The resulting generalized law has a simple form where the effect of memory is clearly parametrized by the real order $\nu$ appearing in \eqref{lo-gen}. For the sake of completness we have fully characterised this model, studying the corresponding relaxation function. The advantage of this generalization is to introduce a physical parameter $\nu$ in the model that allows to consider memory effects in the relaxation process and to find a better agreement with experimental data. Altough the mathematical model is apparently more complicated, it is, in our view, useful to explore the physical applications of fractional Hadamard-like operators with logarithmic kernels. Hence, the main outcomes of our proposal are the following: to consider memory effects in logarithmic-type relaxation models and to start the study on the applications of the Hadamard-like operator \eqref{defO} in the classical equations of mathematical physics. This can be particularly useful in the analysis of diffusive equations related to ultra-slow processes.

   \section{Conclusions}
   
      We have presented  a new approach based on 
      linear integro-differential operators with logarithmic kernel   related to the Hadamard fractional calculus in order to generalize
by a parameter $\nu \in (0,1]$     
    the Lomnitz  creep law.
%% obtained for $\nu=1$.    
For this generalized model   we have derived the constitutive 
stress-strain relation in a form coupling memory effects and time-varying viscosity and also  we have evaluated  numerically  the corresponding relaxation function.
Finally, we note  that the results obtained  in this paper may be useful for fitting  experimental data  in rheology of real materials that exhibit 
   responses in creep and relaxation  varying slower than those of  the Lomnitz model.  
We are thus  confident to have found a suitable application of the Hadamard operators  in the interrelated   fields of   fractional calculus and ultra-slow kinetics.

   \section*{Acknowledgments}
	The work of R. G. and F. M. has been carried out in the framework of the activities of the National Group of Mathematical Physics (GNFM, INdAM).
	The work of G. S. has been partly funded by a research
grant of the Department of Pure and Applied Sciences (DiSPeA)
of the Urbino University “Carlo Bo” (CUP H32I160000000005).	
		The authors would also like to thank  Dr.  Andrea Giusti  for valuable comments and discussion and the anonymous referees for their constructive remarks and suggestions  that helped to improve the manuscript.

%%%%End of the main text
%%%%%%%%%%%%%%%

%%%%%%%APPENDICES
%%%%
\appendix

   	\section{ Hadamard fractional calculus and more general integro-differential operators with logarithmic kernel} \label{app-a}
   
   According to the literature, there are many different definitions of fractional-type derivatives, motivated by the specific mathematical or physical problems under consideration. One of these definitions was
   introduced by Hadamard \cite{had} in 1892 and, although not so frequently used in applications, it is
   generally mentioned in classical reference books on fractional calculus
    (e.g. \cite{samko}, Section 18.3). The Hadamard fractional derivative essentially corresponds to the fractional power of the operator $\delta = \left(t\frac{d}{dt}\right)$ and therefore can be obtained from the well-known Riemann-Liouville fractional derivative, by taking the change of variable $t\rightarrow \ln t$. The Hadamard fractional calculus is now gaining more interest in the mathematical literature (see for example \cite{ram} and references therein), while, as far as we know, few applications in mathematical-physics have been studied. For a survey about this topic, we refer to the textbooks \cite{kilbas} and \cite{samko} and to the review paper of Kilbas \cite{kil}.
 
   We first recall the definition of the Hadamard fractional integral
   \begin{equation}
  (\mathcal{J}^\nu_{a^+}f)(t)= \frac{1}{\Gamma(\nu)}\int_a^t\left(\ln \frac{t}{\tau}\right)^{\nu-1}f(\tau)\frac{d\tau}{\tau}, \quad 0\leq a <t< b\leq \infty.
  \end{equation}
  The left-sided Hadamard fractional derivative of order $\nu\geq 0$ is 
  \begin{equation}
  (\mathcal{D}^\nu_{a^+} f)(t) = \delta^n(\mathcal{J}^{n-\nu}_{a^+}f)(t)=
  \frac{1}{\Gamma(n-\nu)}\left(t\frac{d}{dt}\right)^n\int_a^t\left(\ln \frac{t}{\tau}\right)^{n-\nu-1}f(\tau)\frac{d\tau}{\tau},
  \end{equation}
  where $n = [\nu]+1$ and $[\nu]$ is the integer-part of $\nu$. 
  
  In the recent paper \cite{bal}, the Caputo-like regularization of the Hadamard derivatives was introduced as follows
  \begin{equation}
    (^C\mathcal{D}^\nu_{a^+} f)(t) =(\mathcal{D}^\nu_{a^+})\bigg[f(t)-\sum_{k=0}^{n-1}\frac{\delta^k f(a)}{k!}\left(\ln\left(\frac{t}{a}\right)\right)^k\bigg],
    \end{equation}
   that is equivalent to commute the operator $\delta^n$ with the Hadamard fractional integral (see \cite{bal}, Theorem 1). 
   Therefore, the Caputo-like Hadamard derivative is
   \begin{itemize}
   \item for $\nu\not\in \mathbb{N}$\\
   \begin{equation}%\nonumber
    (^C\mathcal{D}^\nu_{a^+} f)(t)
    = \mathcal{J}^{n-\nu}_{a^+}\left(\delta^n f\right)(t)
    =  \frac{1}{\Gamma(n-\nu)}\int_a^ t
    \left(\ln \frac{t}{\tau}\right)^{n-\nu-1}
    \bigg[\left(\tau\frac{d}{d\tau}\right)^n \! f(\tau)\bigg]\frac{d\tau}{\tau}
    \end{equation}
    \item for $\nu \in \mathbb{N}$\\
    \begin{equation}% \nonumber
     (^C\mathcal{D}^\nu_{a^+} f)(t) = \delta^n f(t) = \left(t\frac{d}{dt}\right)^n f(t).
    \end{equation}
   \end{itemize}
   
   In analogy with this theory, in the recent paper by Beghin et al \cite{jap}, the authors  considered in a probabilistic framework the Caputo-like counterpart of the fractional power of the operator 
   $\widehat{O}^t_1= \left[\left(\frac{a}{b}+t\right)\frac{d}{dt}\right]$.
   Observe that in the special case $a= 0$, we obtain again the Caputo-like regularization of the Hadamard derivative.
   
   In analogy with the Hadamadard derivative, this integro-differential operator can be formally obtained starting from the Caputo derivative by taking the change of variable $t \rightarrow \ln\left(\frac{a}{b}+t\right)$ with the advantage to avoid the singularity for $t=0$. 
   
   The integro-differential operator considered in \cite{jap} is defined:
   \begin{itemize}
   \item if $\nu \not\in \mathbb{N}$ as
   \begin{equation}
   \widehat{O}^t_\nu\,  f(t)= \frac{1}{\Gamma(n-\nu)}\int_{\frac{1-a}{b}}^{t}	\ln^{n-1-\nu}\left(\frac{a+bt}{a+b\tau}\right)
     \widehat{O}^t_n\,  f(\tau)\frac{b}{a+b\tau}d\tau, %\nonumber
   \end{equation}
   \item if $\nu \in \mathbb{N}$ as
   \begin{equation}
   \widehat{O}^t_\nu \,f(t)=  \left(\widehat{O}_1^t \right)^n\, f(t).
   \end{equation}
    \end{itemize}
    
   Some basic results about this operator with logarithmic kernel have been investigated by Beghin et al. in \cite{jap} but a full mathematical analysis about the correct functional setting is still missing. 
   On the other hand, as the authors  have just suggested in that paper, 
   integro-differential operators with logarithmic kernels can have also potential applications in physics.

   We finally observe that Hadamard derivatives and their generalizations can be considered as particular, but interesting cases of fractional derivatives of a function with respect to another function, a topic considered in the book \cite{kilbas} (section 2.5, pp. 99-105). 
   This approach to fractional calculus 
   is poorly known in the mathematical-physics literature; however, quite recently, Almeida \cite{almeida} has studied the regularized Caputo version. 
  
      \section{Numerical approach} \label{app-b}
   Here we present a numerical scheme adopted to solve  the Volterra Integral Equation 
satisfied by our relaxation function $\phi_\nu(t)$ that we re-write henceforth  for the sake of convenience by setting $q=1$, that is  
\begin{equation}
\phi_\nu (t) = 1 - \frac {\nu}{\Gamma(1+\nu)} \int_{0}^{t} \frac{  \ln^{\nu-1}(1+t')     }{1+t'} \phi_\nu(t-t')\, dt'.
\end{equation}

\subsection{Definitions} Our equation  can be rewritten as 

\begin{equation}
\phi_\nu(t) = 1 - \gamma \int_{0}^{t}  K_\nu(t-\tau) \phi_\nu(\tau) d\tau,
\end{equation}
where the Kernel is
\begin{equation}
K_\nu(x) = \frac{\ln^{\nu-1} (1+x)}{1+x} = \frac{1}{\nu} \frac{d}{dx} \ln^\nu (1+x)\,,
\end{equation}
and  the constant $\gamma$ reads 
\begin{equation}
\gamma \equiv  \frac{\nu}{\Gamma(1+\nu)}\,.   
\end{equation}

\subsection{Discretization}

In the following, we  employ the \emph{Euler's product integration}, proposed for the first time by Young in \cite{Young} and described e.g.  in the 1985 monograph by Linz \cite {Linz_BOOK85}.
In particular, we take advantage of the numerical scheme  used
in 1995 by Mainardi et al. \cite {Mainardi-Pironi-Tampieri_1995}
to deal the weakly singular Volterra equation with power law kernel for the  
 generalized Basset problem for a sphere
accelerating in a viscous fluid, see also the 1997 survey by Mainardi
\cite{Mainardi_CISM97}. 
 Here, however, the weak singularity of the kernel is of logarithmic type, 
 that  induces more restrictions on the efficiency of the method.   

The integration interval $(0,t)$ is subdivided into $n$ equal intervals of width $\delta t \equiv h\ll 1$:
\begin{equation}
t_j = jh, \quad j=0, 1, 2, \ldots, n, 
\end{equation}
which allows to write our  Volterra  equation in a piecewice manner:
\begin{equation}
\phi_\nu (t_n) = 1 - \gamma \sum_{j=0}^{n-1} \int_{t_j}^{t_{j+1}} 
\! K_\nu (t_n-\tau)\, \phi_\nu (\tau) \,d\tau.
\end{equation}

We now assume a linear piecewise form for  the function $\phi_\nu(t)$, that is  $\phi_\nu(t) =\phi_j$ for $t_j \le t < t_{j+1}$, which provides
\begin{equation}
\phi_n= 1 - \gamma \sum_{j=0}^{n-1} \phi_j \int_{t_j}^{t_{j+1}}  
\, K_\nu (t_n-\tau)\, d\tau,
\end{equation}
where $\phi_n$ is the approximation for $\phi_\nu(t_n)$.
But
\begin{equation}
\begin{array} {lll}
& {\ds \int_{t_j}^{t_{j+1}} \! \! K_\nu (t_n-\tau) d\tau  
= \int_{t_n-t_{j+1}}^{t_n-t_{j}} \!\! K_\nu (x) \, dx
  =  \int_{ (n-(j+1))h    }^{ (n-j)h}  \! \! K_\nu (x)\, dx}\\ \\
  &= {\ds  \int_{ (n-(j+1))h  }^{ (n-j)h} \! \! \frac{\ln^{\nu-1}(1+x)}{1+x}\,  dx       
   = \frac{1}{\nu} \left. \ln^\nu (1+x)\right  |_{ (n-(j+1))h }^{(n-j)h}}    \\ \\
   &= {\ds \frac{   \ln^\nu ( 1 + (n-j)h )  - \ln^\nu ( 1+ ( n-(j+1))h )} {\nu}}\,.                        
\end{array}
\end{equation}
Hence 
\begin{equation} \label{eq:phin}
\phi_\nu (t_n) = 1 - \gamma \sum_{j=0}^{n-1} \Omega_{n-j} \phi_j\,, 
\end{equation}
with 
\begin{equation}\label{eq:omegan}
 \Omega_{n}  =  \frac{   \ln^\nu ( 1 + nh    )  - \ln^\nu (1+(n-1)h) }{\nu}.
\end{equation}

Equations (\ref{eq:phin}) and (\ref{eq:omegan}) have been implemented in a \texttt{Fortran} program.  The results are depicted in the Figures 1 and 2 in the body of the paper.


\begin{thebibliography}{99}
    
  	  \bibitem{almeida}
		R. Almeida, 
		A Caputo fractional derivative of a function with respect to another function,
		{\it Commun. Nonlinear Sci. Numer. Simulat.} {\bf 44} (2017),
		460--481. 
	
  
      \bibitem{jap}L. Beghin, R. Garra, C. Macci, 
      Correlated fractional counting processes on a finite-time interval. \emph{Journal of Applied Probability} 
      {\bf 52} No 4  (2015), 1045--1061.
      [E-print: {\tt arXiv:1407.6844v2 [math.PR]}] 
      
      \bibitem{Chechkin-et-al_EPL03}
      A. V. Chechkin, J. Klafter,  I. M. Sokolov,
      Fractional Fokker-Planck equation for ultraslow kinetics,
      {\it Europhysics Lett.} {\bf 63} No 3 (2003), 326--332.
      [E-print: {\tt https://arxiv.org/pdf/cond-mat/0301487}] 
      
      \bibitem{WenChen-et-al_FCAA16}
      W. Chen, Y. Liang, X. Hei,
      Structural derivative based on inverse Mittag-Leffler function for modeling ultraslow diffusion,
 {\it Fract. Calc. Appl. Anal.} {\bf   19} No  5 (2016), 1250--1261.
[DOI: 10.1515/fca-2016-0064]      
      
      \bibitem{Christensen_BOOK82}
R.M.  Christensen, 
 {\it Theory of Viscoelasticity}, 2-nd edition,
Academic Press, New York (1982).  [First edition (1972)]
 
      \bibitem{book} 
      R. Gorenflo, A.A. Kilbas, F. Mainardi, S.V. Rogosin, 
      \emph{Mittag-Leffler Functions, Related Topics and Applications}. 
      Springer, Berlin (2014).

\bibitem{Gross_BOOK53}
 B. Gross, %% (1953).
{\it Mathematical Structure of the Theories of Viscoelasticity},
Hermann \& C., Paris (1953).      
      
      
      \bibitem{had} 
      J. Hadamard, 
      Essai sur l’\'etude des fonctions donn\'ees par leur d\'eveloppment de Taylor, 
      \emph{J. Math. Pures Appl.} {\bf 8} (1892), 101--182.
           	 
\bibitem{Hanyga_2005}
A. Hanyga,
Viscous dissipation and completely monotonic relaxation moduli,
{\it Rheologica Acta} {\bf 44} (2005), 614--621.


      
      	\bibitem{bal} F. Jarad, T. Abdeljawad, D. Baleanu,
      	 Caputo-type modification of the Hadamard fractional derivatives, \emph{Advances in Difference Equations} {\bf 2012} No  1 (2012), 
      	 1--8 (2012).
      	 
           \bibitem{Jeffreys_GJRAS58}
H. Jeffreys,
A modification of  Lomnitz's  law of creep in rocks,
{\it Geophys. J. R. Astron. Soc.} {\bf 1} (1958), 92--95.      	 
      	 

      	 
      	 
      	 \bibitem{kil} A.A. Kilbas, 
      	 Hadamard-type fractional calculus,
      	 \emph{J. Korean Math. Soc.}
{\bf   38} No 6 (2001),  1191--1204.
      
     \bibitem{kilbas} 
     A.A. Kilbas, H.M. Srivastava, J.J. Trujillo, 
     \emph{Theory and Applications of Fractional Differential Equations}, Elsevier, Amsterdam (2006).

\bibitem{Linz_BOOK85}
P. Linz, 
{\it Analytical and Numerical Methods for Volterra Equations}
SIAM, Philadelphia  (1985).

      \bibitem{Lo1} C. Lomnitz, Creep measurements in igneous rocks, 
      \emph{J.  Geol.}    {\bf  64} (1956),  473--479.
      
      \bibitem{Lo2} C. Lomnitz, Linear dissipation in solids, 
      \emph{J. Appl, Phys.}      
      {\bf  28} (1957), 201--205.
      
      \bibitem{Lo3} C. Lomnitz, 
      Application of the logarithmic creep law to
      stress wave attenuation in the solid Earth,
       \emph{J. Geophys. Res,.}     {\bf 67} (1962), 365--367.
      
\bibitem{Mainardi_CISM97}
F. Mainardi, 
Fractional calculus, some basic problems in continuum and statistical mechanics,
in :
A. Carpinteri,  and F. Mainardi  (Editors),
{\it Fractals and Fractional Calculus in  Continuum Mechanics},
Springer Verlag, Wien and New York (1997), pp. 291--348.
% (Vol. no 378, series CISM Courses and Lecture Notes, 
% ISBN 3-211-82913-X)
%[Lecture Notes of the Advanced School held at CISM, Udine, Italy,
 %    23-27 September 1996]      
      
      
\bibitem{Mainardi_BOOK10}
F. Mainardi,  \textit{Fractional Calculus and Waves in Linear Viscoelasticity},
Imperial College Press, London (2010). %%, 340 pp.      
      
      \bibitem{Mainardi-et-al_JVC07}
F. Mainardi, A. Mura, R. Gorenflo, M. Stojanovic,
The two forms of fractional relaxation of distributed order,
{\it J. Vibration and Control} {\bf 13} (2007),  1249--1268. %%  (2007).      
[E-print {\tt http://arxiv.org/abs/cond-mat/0701132}]

\bibitem{Mainardi-et-al_JVC08}
F. Mainardi, A. Mura, G. Pagnini, R. Gorenflo,
Time-fractional diffusion of distributed order,
{\it J.  Vibration and Control} {\bf 14} (2008),  1267--1290.      
      [E-print {\tt http://arxiv.org/abs/cond-mat/0701131}]
      
 \bibitem{Mainardi-Pironi-Tampieri_1995}
F. Mainardi, P. Pironi, F. Tampieri.
A numerical approach to the generalized Basset problem for a sphere
accelerating in a viscous fluid, 
in:  P.A. Thibault and D.M. Bergeron (Editors),
{\it Proceedings CFD 95}, Vol. II, pp. 105-112 (1995).
[3-rd Annual Conference of the Computational Fluid Dynamics Society
of Canada, Banff, Alberta, Canada, 25-27 June 1995].      
      
      \bibitem{Main} F. Mainardi, G. Spada, 
  On the viscoelastic characterization    of the Jeffreys–-Lomnitz law of creep, 
            \emph{Rheol. Acta} {\bf 51} (2012), 783--791.
            [E-print: {\tt http://arxiv.org/abs/1112.5543}]
            
  \bibitem{geo} G. Mavko, N. Saxena, 
  Rock-physics models for heterogeneous creeping rocks and viscous fluids, {\it Geophysics} {\bf 81(4)} (2016), D427-–D440.  
      
  \bibitem{Metzler-Klafter_JPhysics04}
 R. Metzler,  J.  Klafter,  %%%  (2004).
The restaurant at the end of the random walk: Recent developments
 in the description of anomalous transport by fractional dynamics,
 {\it J. Phys. A. Math. Gen.}  {\bf 37} (2004),  R161--R208.      
      
\bibitem{Molinari_1973}
 A. Molinari,  
Visco\'elasticit\'e lin\'eaire et fonctions compl\`etement monotones,
{\it Journal de M\'ecanique} {\bf 12} (1973),  541--553.
      
      
       \bibitem{Holm} V. Pandey, S.Holm, 
       Linking the fractional derivative and the Lomnitz creep law
            to non--Newtonian time--varying viscosity, 
            \emph{Physical Review E}     {\bf 94} (2016), 032606/1-6.
            
 \bibitem{Pipkin_BOOK86}
 A.C. Pipkin, % (1986).
{\it Lectures on Viscoelastic Theory},
2-nd Edition,
Springer Verlag, New York (1986).  [First Edition 1972]
%% (Applied Mathematical Sciences No 7)

		\bibitem{samko} S.G. Samko, A.A. Kilbas, O.I. Marichev, \emph{Fractional Integrals and Derivatives: Theory and Applications}, Gordon and Breach (1993).
            
     \bibitem{ram} R.K. Saxena, R. Garra, E. Orsingher,
      Analytical solution of space-time fractional telegraph-type equations involving Hilfer and Hadamard derivatives,
      \emph{Integral Transforms and Special Functions} 
      {\bf 27} No 1 (2016),  30--42.
                  
              \bibitem{Sinai_TBA83}
              Ya. G. Sinai,
  The limiting behavior of a one-dimensional random walk
                 in a random medium,
{\it Theory Probab. Appl.} {\bf 27} No 2 (1983), 256--268.  
 
            \bibitem{Scheidegger_1970a}
            A.E Scheidegger,
            On the rheology of rock creep,
            {\it Rock Mechanics} {\bf 2} (1970), 138--145.
            
                      \bibitem{Scheidegger_1970b}
            A.E Scheidegger,
           The rheology of the Earth in the intermediate time range, 
            {\it Annals of Geophysics} {\bf 23} No 1 (1970),  27--43.
            
       \bibitem{sch} R.L. Schilling, R. Song, Z. Vondracek,
        \emph{Bernstein Functions: Theory and Applications}, 
        %Vol. 37. Walter 
        De Gruyter (2012).
       
       \bibitem{Strick_JGR84}
E. Strick,
Implications of Jeffreys--Lomnitz transient creep,
{\it J. Geophys. Res.} {\bf 89} (1984), 437--451.

      \bibitem{Young} A.Young, 
      Approximate product-integration,
      {\it Proc.R.Soc. Lond. Ser.A}, {\bf 224} (1954), 552--561.

\end{thebibliography}
        \end{document}